\documentclass[]{interact}

\usepackage{subfigure}% Support for multi-panel figures

\usepackage[numbers,sort&compress]{natbib}% Numerical citation support
\bibpunct[, ]{[}{]}{,}{n}{,}{,}
\renewcommand\bibfont{\fontsize{10}{12}\selectfont}

\begin{document}

\title{The Geometry of Wigwag and Reversal Linework in Blown Glass}

\author{
\name{Anduriel Widmark\thanks{CONTACT Anduriel Widmark. Email: anduriel@andurielstudios.com}}
\affil{Independent researcher, Denver, Colorado, USA}
}

\maketitle

\begin{abstract}
Wigwag and reversal linework in blown glass produces periodic wave, spiral, and radiating patterns by twisting longitudinal stripes in alternating directions and then reorganizing the patterned tube. Although glassworkers have developed many consistent techniques and identifiable pattern forms, these patterns do not appear to have been widely described within a common coordinate framework.

This paper presents one way of modeling these patterns mathematically. Periodic linework is defined on a cylinder, represented on a sphere, and read relative to a switched axis. The switched pattern may then be re-formed into cylindrical or circular configurations and further manipulated into square, rectangular, or compound forms. The circular case provides a particularly simple radial representation of the switched spherical field. This sequence describes the characteristic organization of a reversal, with the same geometry extending to stacked encalmo sections and other compound structures. The goal is to provide a concise mathematical language for interpreting this broad family of patterned glass forms.
\end{abstract}


\section{Introduction}

Wigwag linework is a family of patterned-glass techniques in which longitudinal stripes are repeatedly twisted in alternating directions, producing periodic wigs and wags. Further manipulation, through condensation, joining, reorientation, or reversal, can produce radial arrangements, spiral-like terminations, and compound patterns that appear entirely different from the initial cylindrical linework.

The glasswork is technically demanding and can require a high level of refinement. The glassworker must control rotation, heat distribution, wall thickness, line spacing, and timing while manipulating molten glass. Small variations can propagate through many repeated lines, requiring considerable skill to maintain precise symmetry and registration.

Contemporary glassworking literature provides an important practical record of these techniques, including longer-standing practices and more recent developments within the borosilicate glass art movement. Bandhu Scott Dunham's
\textit{Contemporary Lampworking} documents wigwags and a broad vocabulary of patterned-glass processes within modern lampworking \cite{Dunham2025,Dunham2010}. Such patterns are commonly described through process: line count, direction and amount of twist, construction of lined tubing, joining of sections, changing the axis, and recombining patterned material.

Here, a geometric description asks how the symmetries of the linework are maintained and transformed as it is reorganized in space.

This work presents one way of modeling several characteristic wigwag and reversal structures within a shared coordinate framework. Photographs of physical glass examples accompany the diagrams so that each mathematical stage can be compared directly with realized objects.

The starting point is an evenly distributed set of longitudinal lines on a cylinder undergoing periodic angular displacement. Periodic examples are used here because they make symmetry and coordinate relationships particularly clear; aperiodic wigwags are also used in practice and are not excluded by the framework. A reversal is then described by carrying the pattern onto a sphere, changing its reference axis, and re-forming the switched field into cylindrical, circular, or other patterned configurations.

The examples presented are used as symmetric reference cases, not as ideal or preferred physical forms, but as convenient choices where periodicity, symmetry, and axis switching are particularly clear. Actual glass allows much greater variation through changes in heating, diameter, compression, joining, and deliberate manipulation.

\section{Periodic Linework on a Cylinder}

Let the source tube have angular coordinate
\(\theta\in[0,2\pi)\) and normalized longitudinal coordinate
\(v\in[0,1]\). For \(N\) evenly distributed lines, a convenient symmetric
wigwag is

\begin{equation}
\theta_k(v)
=
\frac{2\pi k}{N}
+
A\sin\left(
2\pi f\left(v-\frac{1}{2}\right)+\phi
\right),
\qquad
k=0,\ldots,N-1.
\label{eq:wigwag}
\end{equation}

Here \(N\) is line count, \(A\) is angular amplitude, \(f\) is longitudinal
frequency, and \(\phi\) is phase.

Equation~\eqref{eq:wigwag} describes the organization of the resulting
linework. Other periodic functions could represent sharper or less regular
changes of direction. The sine form simply isolates the main variables.

Because the source is cylindrical, \(\theta\equiv\theta+2\pi\). The left and
right boundaries of an unwrapped diagram represent the same location on the
physical tube, so line identity remains cyclic.

Changing \(N\) changes line count, \(A\) controls angular travel around the
cylinder, \(f\) determines the number of oscillations along its length, and
\(\phi\) shifts the position of the periodic pattern along the length of the
tube.

\begin{figure}[htbp]
\centering

\subfigure[]{%
\includegraphics[width=0.47\textwidth]
{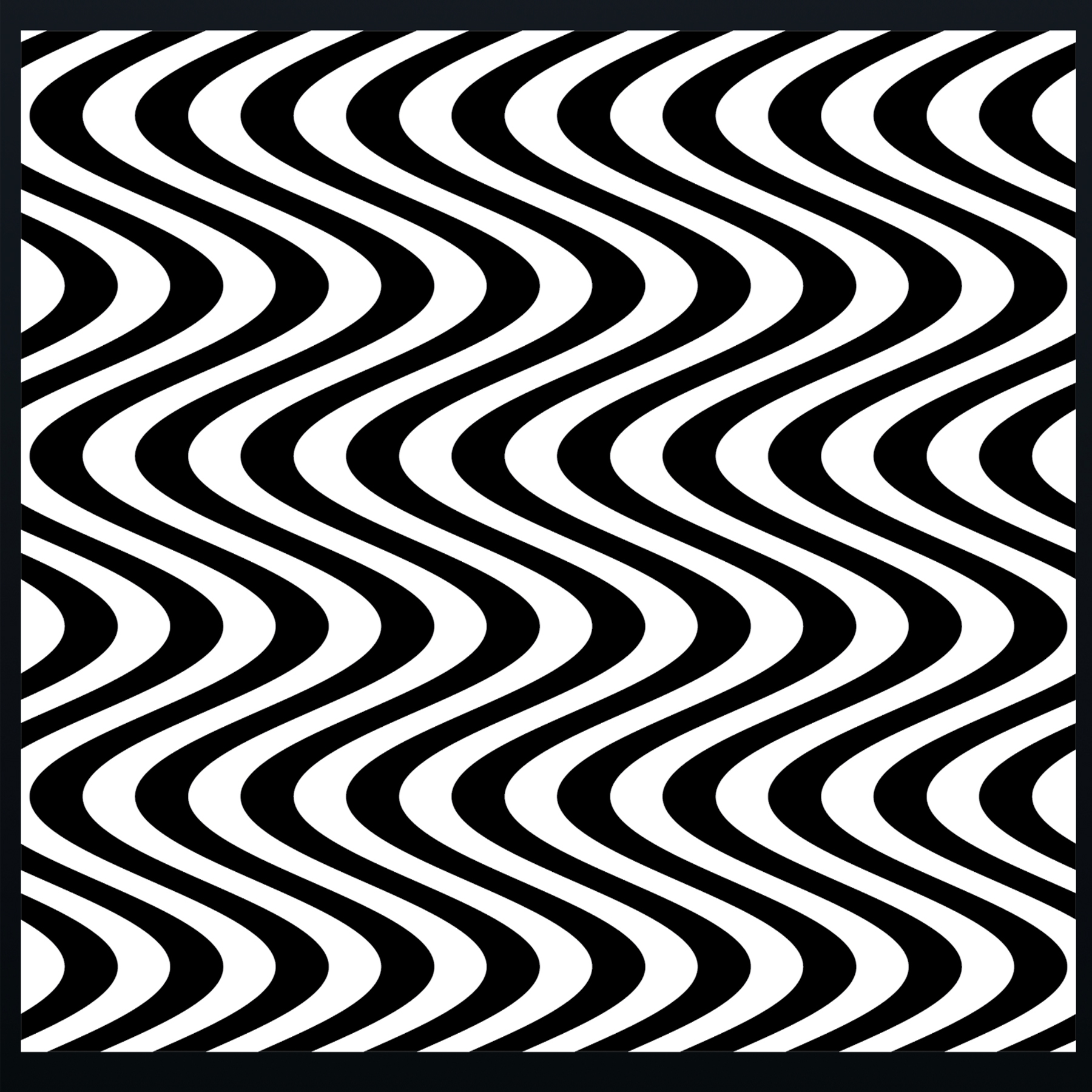}
\label{fig:1a}
}
\hfill
\subfigure[]{%
\includegraphics[width=0.47\textwidth]
{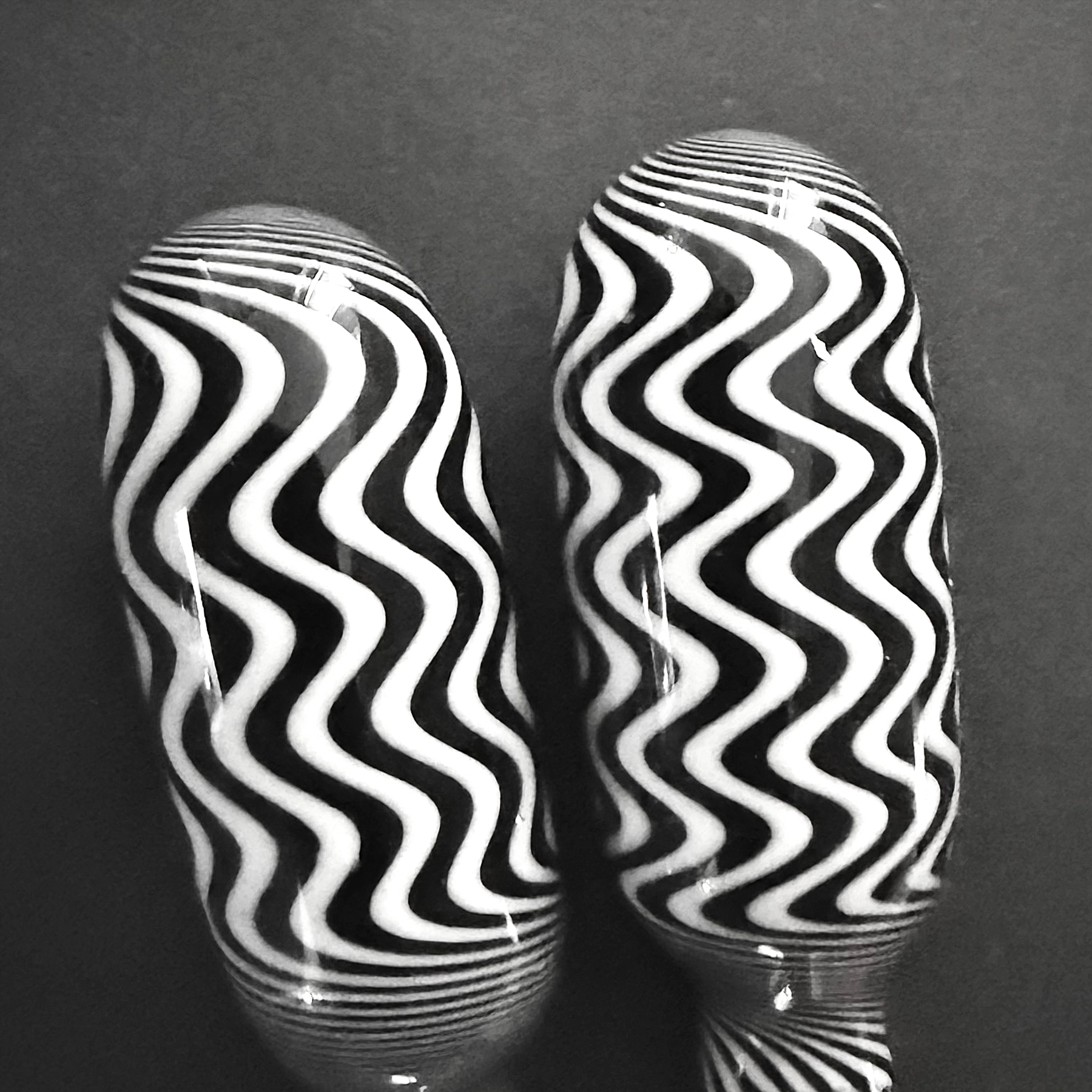}
\label{fig:1b}
}

\caption{(a) Digital unwrapped wigwag representation; 
(b) author-made wigwag glass cylinders.}

\label{fig:wigwag-cylinder}
\end{figure}

\section{The Cylinder Becomes a Sphere}

The reversal becomes clear after transferring the cylindrical coordinates to a
sphere. With \(z=2v-1\),

\begin{equation}
P(\theta,z)
=
\left(
\sqrt{1-z^2}\cos\theta,\,
\sqrt{1-z^2}\sin\theta,\,
z
\right).
\label{eq:sphere}
\end{equation}

The angular coordinate is retained while longitudinal position becomes the
vertical spherical coordinate.

Near \(z=0\), the pattern lies around the equatorial region of the sphere.
Toward either end, the circular cross-section contracts. At \(z=\pm1\), an
entire end of the cylindrical coordinate domain converges to a point.

This gives a geometric origin for the termination. Lines formerly separated
around the circumference of the tube are brought together at a spherical pole.

The specific appearance of that region remains open to artistic control. A
glassblower may tighten, loosen, rotate, compress, or otherwise reshape the
innermost spiral. Here, termination character is considered an additional
artistic parameter.

\begin{figure}[htbp]
\centering

\subfigure[]{%
\includegraphics[width=0.47\textwidth]
{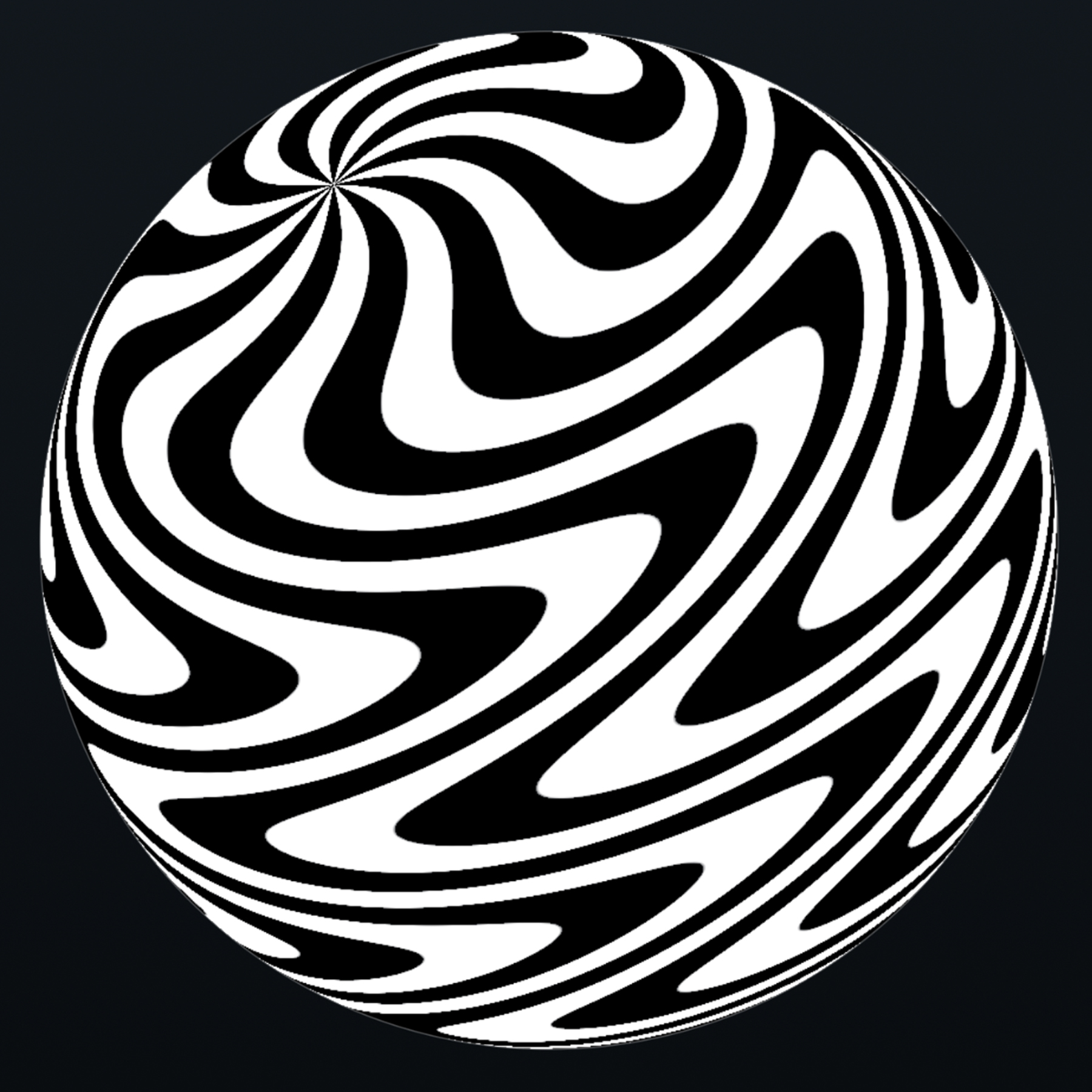}
\label{fig:2a}
}
\hfill
\subfigure[]{%
\includegraphics[width=0.47\textwidth]
{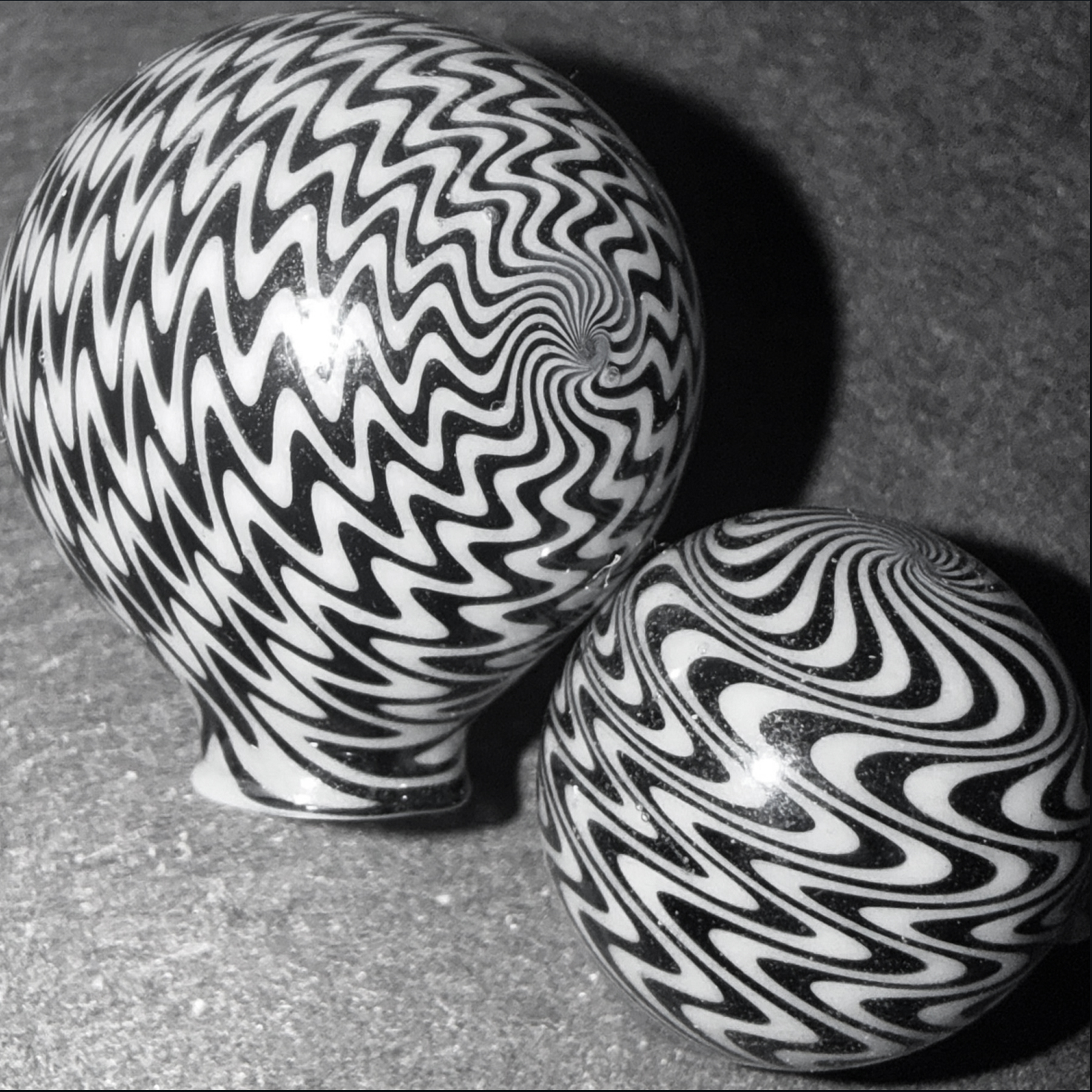}
\label{fig:2b}
}

\caption{(a) Digital spherical wigwag; 
(b) two author-made spherical wigwag forms.}

\label{fig:wigwag-sphere}
\end{figure}

\section{Axis Switching and Forming the Reversal}

The spherical stage by itself does not produce the characteristic reversal.
The essential operation in this model is a change of reference axis.

Let

\[
n=(\sin\alpha,0,\cos\alpha)
\]

be the switched axis, with perpendicular basis vectors

\[
e_1=(\cos\alpha,0,-\sin\alpha),
\qquad
e_2=(0,1,0).
\]

Coordinates relative to this new orientation are

\begin{equation}
u=P\cdot n,
\qquad
\psi=
\operatorname{atan2}
\left(
P\cdot e_2,\,
P\cdot e_1
\right)
-\psi_0.
\label{eq:axis-switch}
\end{equation}

The angle \(\alpha\) gives the change between the original and switched axes,
while \(\psi_0\) selects the angular seam.

This is the central distinction between reversing the direction of a wigwag and
producing an axis-switch reversal. Changing the sign of the original
oscillation changes its direction on the same cylinder.
Equation~\eqref{eq:axis-switch} instead changes the orientation from which the
complete spherical pattern is read.

When \(\alpha=0\), the two axes coincide. When \(\alpha=\pi/2\), they are
perpendicular. Intermediate values describe oblique switches.

After the axis change, the same reorganized spherical field may be realized in
several physical forms. It may be re-formed around a new cylindrical axis,
compressed into a circular field, or manipulated further into square,
rectangular, and compound configurations. The switched spherical field
provides a geometric intermediary between the source construction and these
realized forms, the exact final form is determined by the forming of the glass.

A particularly useful symmetric planar representation is obtained by
expressing the switched field radially. Let

\begin{equation}
r=\sqrt{\frac{1+u}{2}},
\qquad
(p,q)=
\left(
r\cos\psi,\,
r\sin\psi
\right).
\label{eq:radial}
\end{equation}

Using \(1-u\) in place of \(1+u\) chooses the opposite switched pole as the
center.

Equation~\eqref{eq:radial} gives a normalized azimuthal equal-area
representation of the switched spherical field. It is the Lambert azimuthal
equal-area construction scaled to a unit disk; the normalization changes the
overall scale but not the relative-area property \cite{Snyder1987}. This
circular field provides a simple mathematical comparison with the physical
reversal-window form.

In the perpendicular reference case, \(\alpha=\pi/2\), the two termination
regions inherited from the original tube satisfy \(u=0\). They occur at the
same radius in the circular field and at opposite angular positions.

The circular representation does not uniquely determine the physical geometry
of a finished reversal. A switched field may instead be re-formed
cylindrically or deformed further as the glass is compressed, stretched,
joined, or reshaped.

A smooth disk-to-square mapping provides a useful symmetric reference for
square and rectangular reversals. Here an elliptical-grid transformation is
used because it preserves the broad radial organization of the circular field
while smoothly redistributing it toward a four-sided boundary. For disk
coordinates \((p,q)\) and normalized square coordinates \((X,Y)\), a compact
form of the correspondence is

\begin{equation}
p=X\sqrt{1-\frac{Y^2}{2}},
\qquad
q=Y\sqrt{1-\frac{X^2}{2}}.
\label{eq:elliptical-grid}
\end{equation}

Equation~\eqref{eq:elliptical-grid} is written in the square-to-disk direction;
its inverse carries the circular field into the square. For a rectangular
form, the resulting \(X\)- and \(Y\)-coordinates may be rescaled independently
after the square transformation \cite{Lambers2016}.

As with the preceding representations, this mapping describes the geometric
organization of a symmetric reference case rather than the exact mechanical
deformation of every physical object. Related smooth disk-to-square mappings
can produce similar overall pattern organizations with different local
distributions of distortion.

The essential geometry of the model is

\[
\text{periodic cylinder}
\longrightarrow
\text{sphere}
\longrightarrow
\text{axis switch}
\longrightarrow
\text{physical re-forming and further deformation}.
\]

The resulting form may be cylindrical, circular, square, rectangular,
compound, or shaped further.

\begin{figure}[htbp]
\centering

\subfigure[]{%
\includegraphics[width=0.47\textwidth]
{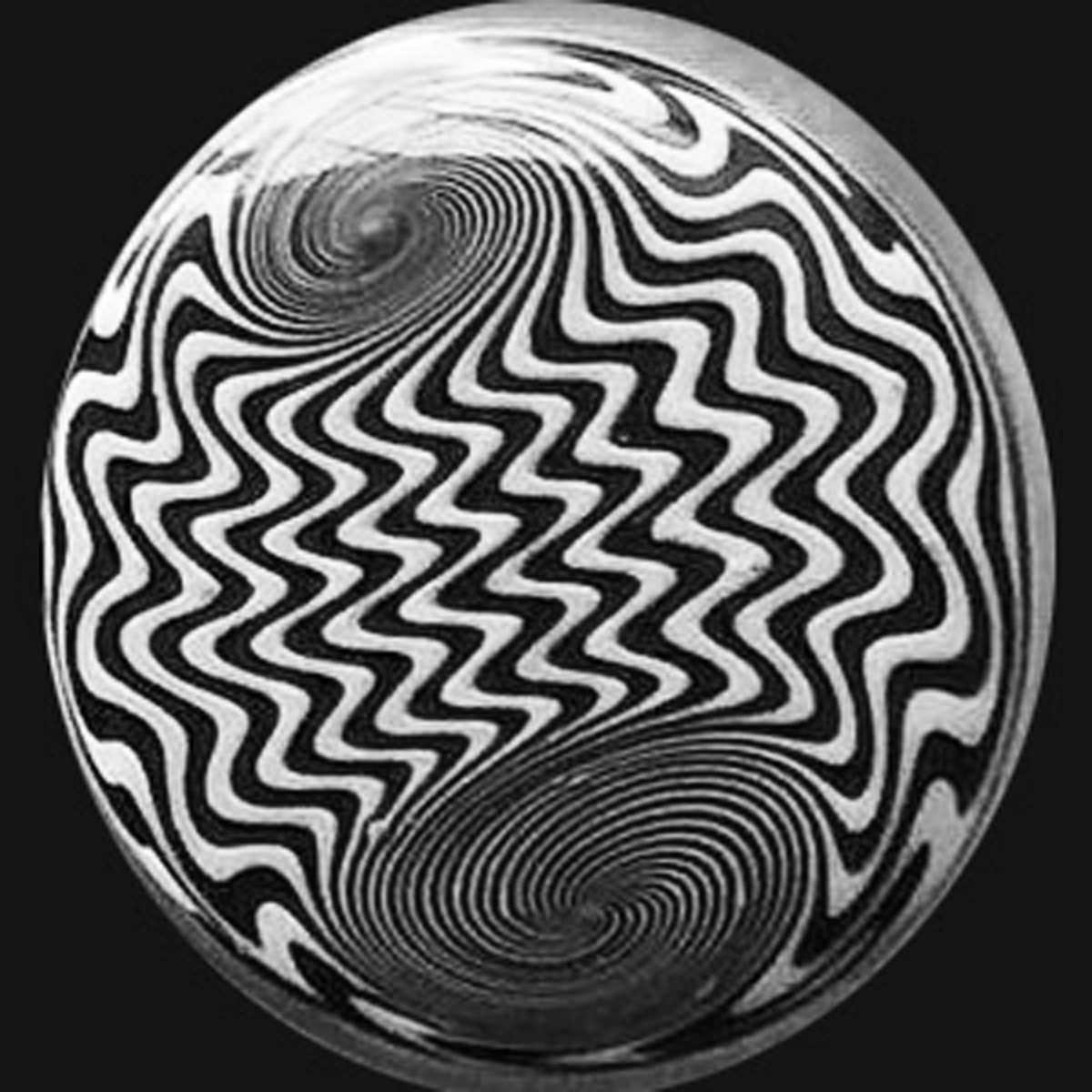}
\label{fig:3a}
}
\hfill
\subfigure[]{%
\includegraphics[width=0.47\textwidth]
{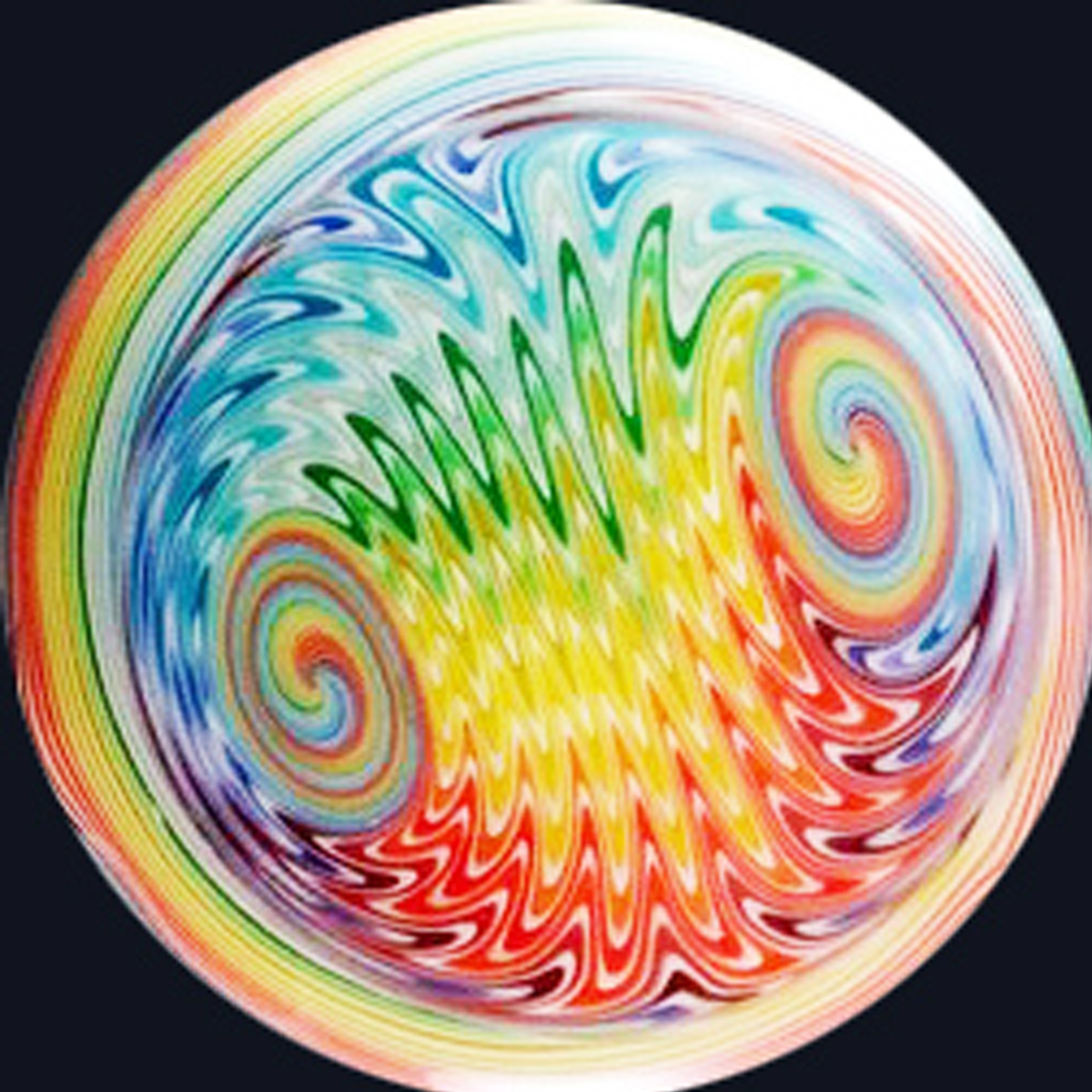}
\label{fig:3b}
}

\caption{(a) and (b) two different author-made circularly compressed reversals/ windows.}

\label{fig:wigwag-glass}
\end{figure}

\begin{figure}[htbp]
\centering

\subfigure[]{%
\includegraphics[width=0.47\textwidth]
{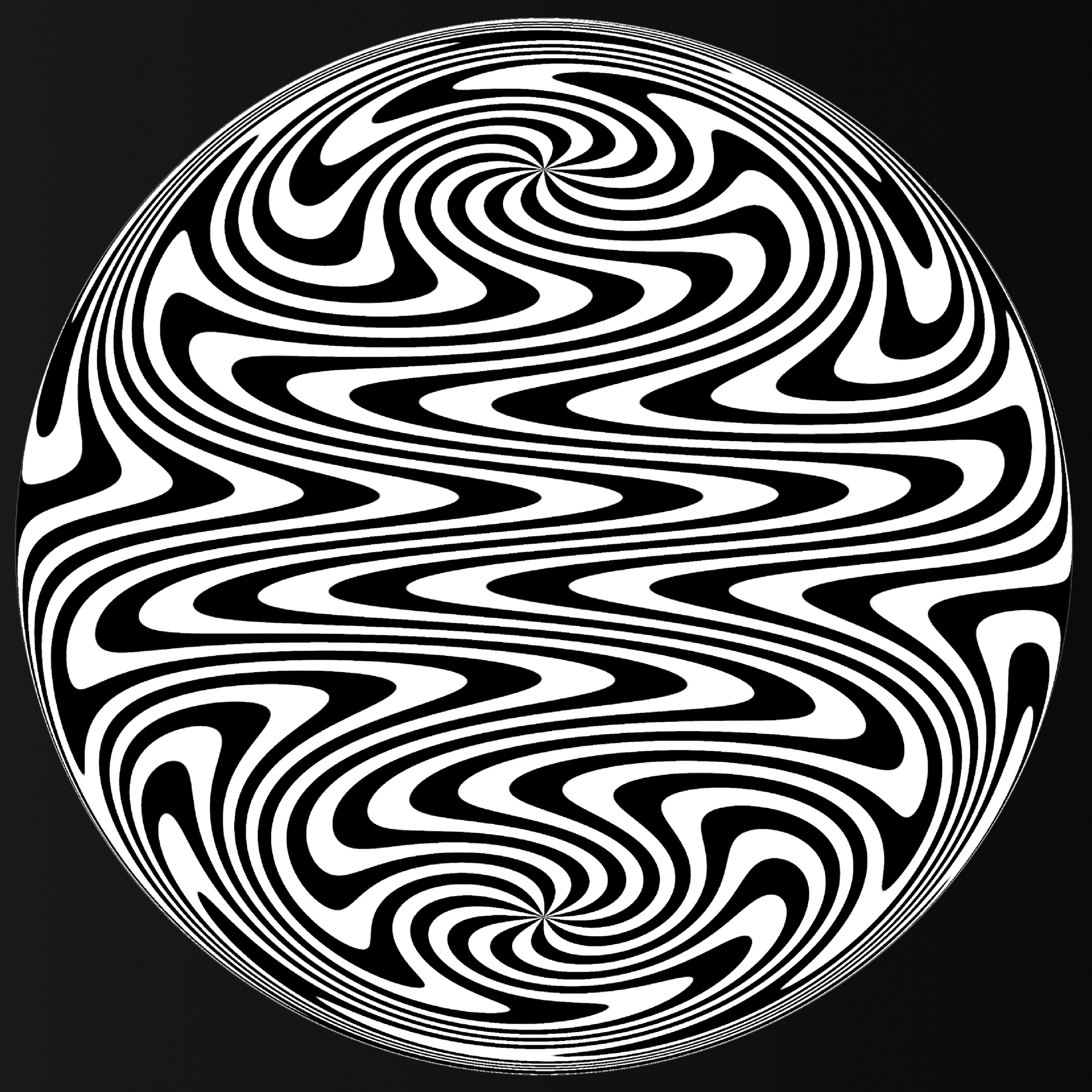}
\label{fig:4a}
}
\hfill
\subfigure[]{%
\includegraphics[width=0.47\textwidth]
{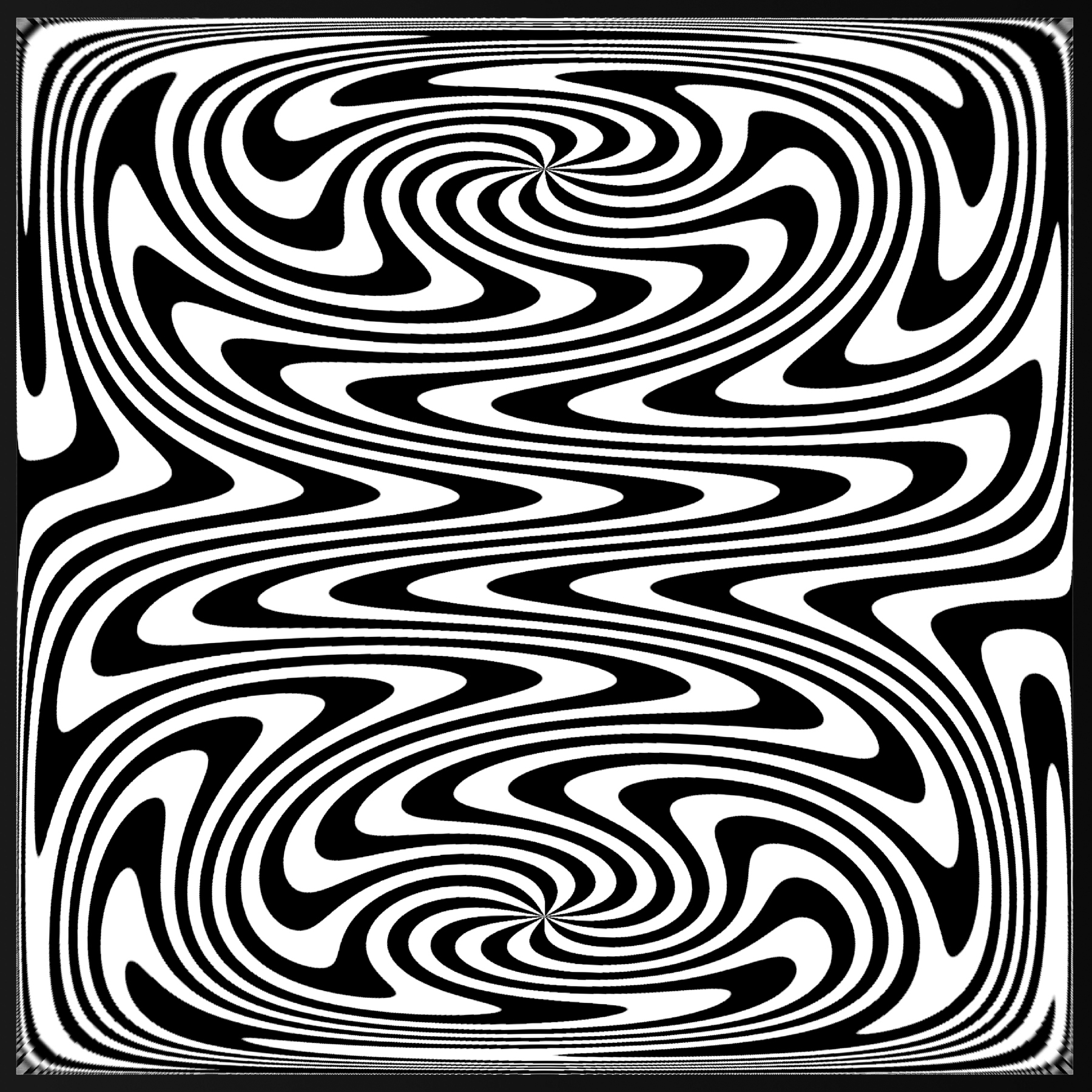}
\label{fig:4b}
}

\caption{(a) circular digital mapping; (b) square mapping of the same reversal pattern.}

\label{fig:wigwag-reversal}
\end{figure}

\section{Encalmo Stacking and Compound Constructions}

A patterned source does not need to remain uniform along its length. Distinct
regions may be joined as encalmo sections, allowing portions of one composite
tube to carry different line counts, amplitudes, frequencies, phases, colors,
or twist directions.

Partition \(v\) into intervals \(I_j=[a_j,a_{j+1}]\). Within section \(j\),

\begin{equation}
t_j
=
\frac{v-a_j}{a_{j+1}-a_j},
\qquad
W_j(t_j)
=
A_j\sin\left(
2\pi f_j\left(t_j-\frac{1}{2}\right)+\phi_j
\right)+c_j.
\label{eq:encalmo}
\end{equation}

Here \(t_j\) normalizes the local section to \([0,1]\), so \(\phi_j\) gives
phase relative to the midpoint of that section.

The offset \(c_j\) may be chosen to make the angular displacement continuous
across the joint, while each region may also have its own line count \(N_j\).

Two possible cases are:

\textbf{Continuous-field reversal.}
The assembled encalmo tube is treated as one source and undergoes a single
spherical mapping and axis switch. Each section retains its distinct pattern
parameters within the unified transformation.

\textbf{Section-wise reversal and reassembly.}
Each encalmo section is treated as a local source and reversed independently
before the resulting sections are reassembled. The reassembled tube may then
serve as a new source for a further shared sphere-and-axis reversal or other
subsequent forming.

Once segmentation and recombination are available, the possible structures
grow rapidly. Sections may vary in line count, periodicity, color, direction,
phase, or orientation. A joined pattern can also be re-wagged, meaning
subjected to another alternating twist, so that one periodic organization is
carried through a second periodic operation. Reversal, re-wagging, further
sectioning, and subsequent joining can be applied in many orders.

The apparent complexity of a finished object arises from repeated
transformations. A relatively small vocabulary of operations, including
periodic twisting, segmentation, axis switching, re-wagging, forming, and
recombination, can be composed into an open-ended family of finite compound
structures.

In practice, borosilicate tubing lends itself particularly well to this kind of
geometric composition because the surface design of patterned tubing can be
cut, joined, rotated, stretched, reoriented, compressed, and manipulated
through multiple stages before being incorporated into a continuous final form
\cite{Dunham2025,Dunham2010}.

\begin{figure}[htbp]
\centering

\subfigure[]{%
\includegraphics[width=0.47\textwidth]
{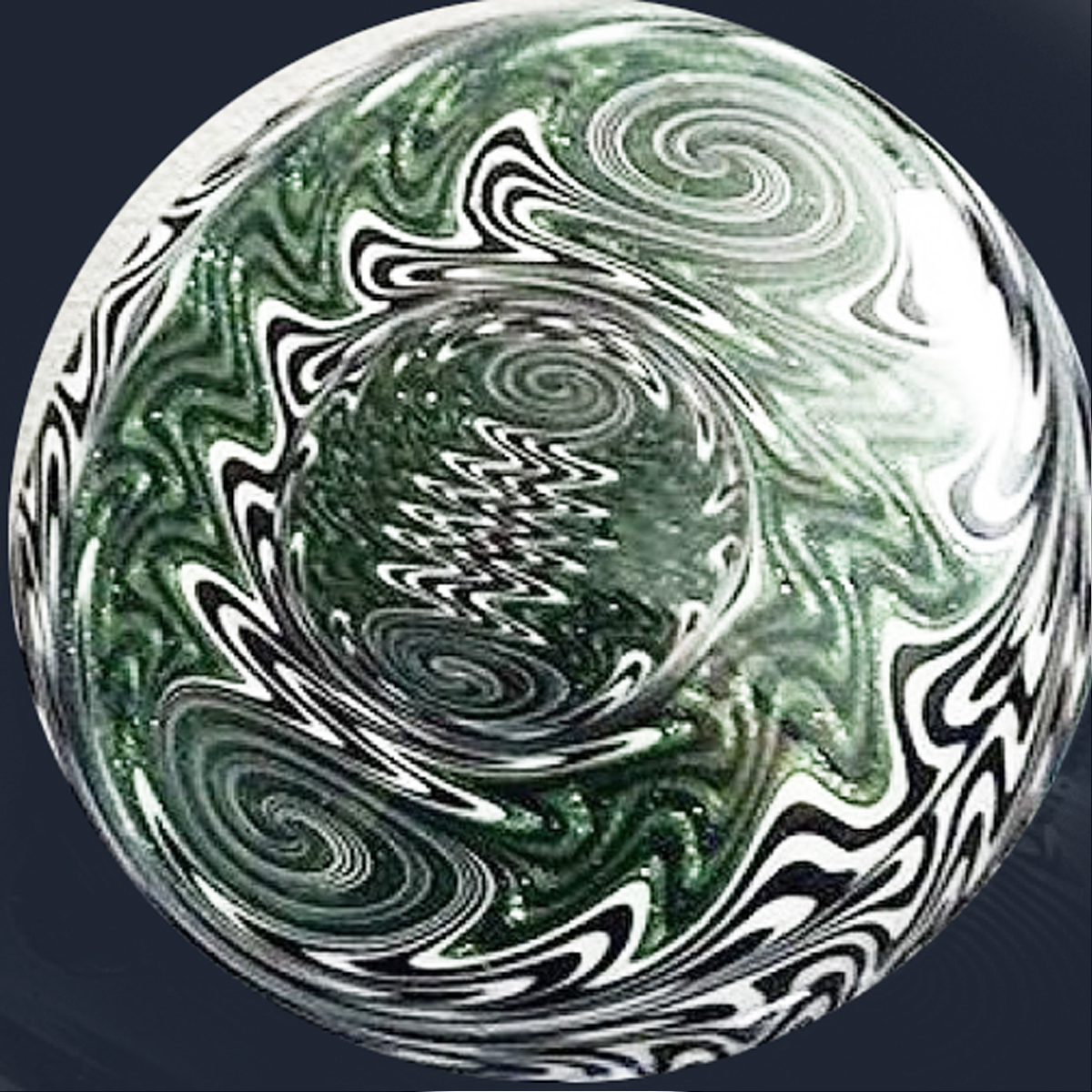}
\label{fig:5a}
}
\hfill
\subfigure[]{%
\includegraphics[width=0.47\textwidth]
{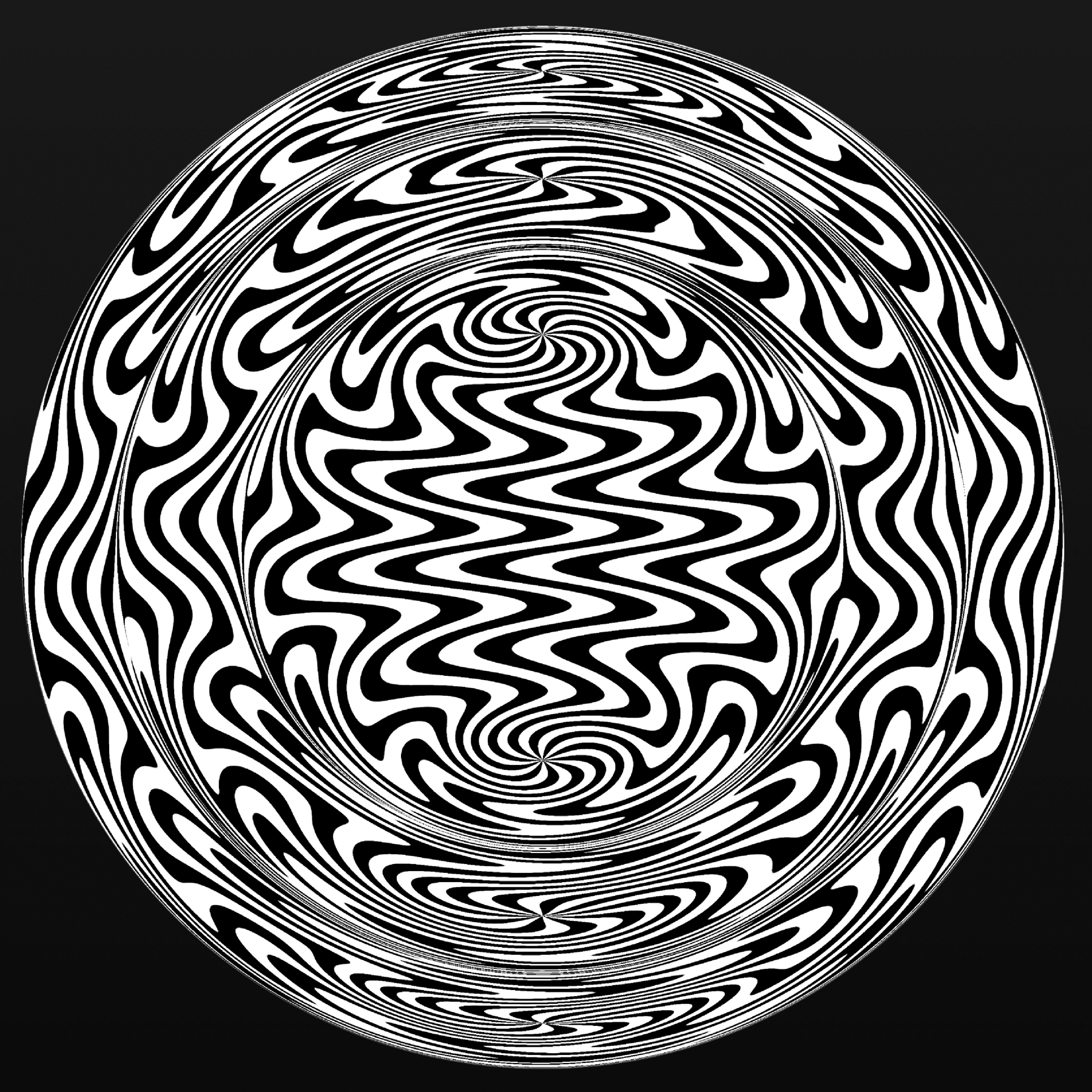}
\label{fig:5b}
}

\caption{(a) author-made two-section encalmo pendent; (b) digital three-section circularly compressed encalmo model.}

\label{fig:wigwag-encalmo}
\end{figure}

\section{Relation to Conformal and Inversive Geometry}

The radial and spiral organization of reversal linework invites comparison with
conformal and inversive geometry. Stereographic projection and M\"obius
transformations provide classical conformal and inversive comparisons, while
rational maps and elliptic functions belong to broader complex-analytic
families capable of producing related spherical transformations.

The coordinate construction described here is not conformal as a whole. The
change of spherical reference axis is a rotation and preserves the geometry of
the sphere. The radial representation in Equation~\eqref{eq:radial} is
equal-area rather than angle-preserving, while the elliptical-grid
correspondence in Equation~\eqref{eq:elliptical-grid} is neither conformal nor
equal-area \cite{Lambers2016}.

These mathematical frameworks provide comparisons and possible extensions
rather than explanations of the basic reversal. Once linework has been placed
on a sphere, conformal, inversive, and other transformations can offer
additional ways to reorganize or display it.

Schleimer and Segerman demonstrate the use of M\"obius transformations,
rational maps, elliptic functions, and related methods for transforming
spherical imagery \cite{SchleimerSegerman2016}. Earlier work by Oscar S. Adams
used elliptic functions in conformal world-map construction \cite{Adams1925},
a related mathematical lineage also discussed by Schleimer and Segerman.
These examples place the present geometry within a broader family of spherical
transformations while preserving the distinction between the equal-area radial
representation, the additional disk-to-square deformation used here, and
conformal mappings.

\section{Scope of the Geometric Description}

The equations describe spatial organization rather than the full mechanics of
molten glass. Temperature, viscosity, wall thickness, diameter, compression,
inflation, pulling, joining, and deliberate asymmetry can deform the symmetric
reference cases while leaving their geometric organization recognizable.

The mappings used here describe relationships among surface coordinates rather
than a strain-free mechanical unfolding of molten glass. Subsequent forming
may introduce additional stretching, compression, local deformation, and
changes in wall thickness while preserving recognizable organization of the
switched field.

More detailed termination behavior, curved vessel surfaces, transparent or
layered linework, interacting periodic patterns, and their optical effects
provide many possible directions for further investigation.

\section{Conclusion}

Wigwag and reversal linework can be described through a compact sequence of
geometric operations. A periodic pattern is established on a cylinder.
Transferring it to a sphere brings its end boundaries toward poles, and
changing the spherical reference axis reorganizes the linework into a new
orientation.

The switched spherical field provides a geometric organization from which
cylindrical, circular, square, rectangular, and compound reversal forms can be
developed through subsequent forming and recombination. The normalized radial
representation provides a particularly simple symmetric description of the
circular, or window, case, while the elliptical-grid correspondence gives a
compact reference construction for carrying that organization into a square or
rectangular field.

This description separates geometric structure from artistic and material
choice. Line count, color, frequency, amplitude, termination character,
asymmetry, boundary shape, and physical deformation can change dramatically
while the underlying coordinate relationships remain recognizable.

Encalmo stacking extends the geometry into compound structures. Patterned
sections can be joined, reversed together or independently, re-wagged, formed
into different boundaries, and recombined. Visually complex results can be
understood as developing through repeated use of a relatively small set of
spatial operations.

Taken together, these constructions provide a geometric basis for describing a
broad family of wigwag and reversal patterns within a shared coordinate
framework. Patterns that glassworkers distinguish through process, appearance,
or varied terminology can also be compared through their underlying
relationships of frequency, amplitude, axis, surface, boundary, and
transformation. The framework complements the practical and material language
of glassworking with a precise description of spatial structures that can
persist through substantial variation in form and process. Wigwag and reversal linework can be understood both as a highly variable
material practice and as a coherent family of geometric transformations.

\section*{Disclosure statement}

the author did not report competing interests.

\section*{Funding}

This research did not receive external funding.

\section*{Declaration of generative AI use}

The author used ChatGPT (GPT-5.6 Sol) for language editing and text formatting.

\bibliographystyle{tfs}
\bibliography{references}

\end{document}